\input amstex

\documentstyle{amsppt}
\NoRunningHeads

\topmatter

\title
 Almost continuous orbit equivalence for  non-singular
 homeomorphisms
 \endtitle
\author  Alexandre~I.~Danilenko and Andr\'es del Junco
\endauthor

\address
 Institute for Low Temperature Physics
\& Engineering of National Academy of Sciences of Ukraine, 47 Lenin
Ave.,
 Kharkov, 61164, UKRAINE
\endaddress
\email danilenko\@ilt.kharkov.ua
\endemail

\address
Department of Mathematics, University of Toronto, M5S 3G3, CANADA
\endaddress

\email deljunco\@math.toronto.edu
\endemail

\thanks The first named author thanks the University of Toronto for
supporting his visit to Toronto,
where the main part of this work was done.
\endthanks

 \abstract
Let $X$ and $Y$ be Polish spaces with non-atomic Borel measures $\mu$
and $\nu$ of full support.
Suppose that $T$ and $S$ are ergodic non-singular homeomorphisms of
$(X,\mu)$ and $(Y,\nu)$ with
continuous Radon-Nikodym derivatives. Suppose that either they are
both  of type $III_1$ or that they are both of type $III_\lambda$,
$0<\lambda<1$ and, in the  $III_\lambda$ case, suppose in addition
that both `topological  asymptotic ranges' (defined in the article)
are $\log\lambda\cdot\Bbb Z$. Then there exist
invariant dense
$G_\delta$-subsets  $X'\subset X$ and $Y'\subset Y$ of full measure
and a non-singular homeomorphism $\phi: X' \to Y'$ which is an orbit
equivalence between $T|_{X'}$ and $S|_{Y'}$,   that is
$\phi\{T^{i}x\} = \{S^{i}x\}$ for all $x \in X'$. Moreover
  the Radon-Nikodym derivative
$d\nu\circ\phi/d\mu$ is continuous on $X'$ and, letting $S' =
\phi^{-1}S \phi$ we have $Tx= {S'}^{n(x)}x$ and $S' = T^{m(x)}x$ where
$n$ and $m$ are continuous on $X'$.
  \endabstract

 \endtopmatter
   \document

\head 0. Introduction
\endhead

The study of {\it measure-theoretic orbit equivalence} was
initiated by H.~Dye \cite{Dy1}, \cite{Dy2} who proved that the orbit
equivalence relations of any two ergodic finite measure preserving
transformations are  measure theoretically isomorphic  when
restricted to some invariant subsets of full measure. W.~Krieger
obtained  in \cite{Kr1}, \cite{Kr2} a far-reaching generalization of
Dye's theorem and classified, up to measure-theoretic orbit
equivalence, all ergodic non-singular transformations.  Dye's and
Krieger's theorems extend to  ergodic actions of countable
amenable groups \cite{CFW}. In contrast to that, any non-amenable
 countable group has uncountably many non-orbit equivalent
probability preserving ergodic actions \cite{Hj}, \cite{Ep}. On the
other hand, amenability plays no role in {\it generic orbit
equivalence}: given any two topologically transitive countable groups
of homeomorphisms of perfect Polish spaces, their orbit equivalence
relations are  homeomorphic after restriction to  some invariant
dense $G_\delta$ subsets \cite{SWWr}.

Consider now the class of    homeomorphisms
of Polish spaces equipped  with ergodic nonsingular measures of full
support
such that the Jacobian (Radon-Nikodym derivative) of the
homeomorphism with respect to its non-singular measure is continuous. We consider a
notion of orbit equivalence for these which combines the notions of generic and
measure-theoretic orbit equivalence. We say
that two such homeomorphisms are {\it almost continuously orbit equivalent} if
their orbit equivalence relations restricted to some invariant dense
$G_\delta$-subsets of full measure are isomorphic by a nonsingular
homeomorphism which has a continuous Radon-Nikodym derivative and
such that the co-cycles associated to the orbit equivalence are
continuous  (see
Definition~1.1 below).

T.~Hamachi  and M.~Keane initiated
investigation of such a kind of equivalence in the measure-preserving
case
\cite{HK} by constructing  a weakly  continuous (finitary in
their terminology) orbit
equivalence between the binary  and ternary  odometers. We use the
term `weakly' here since they  only required that the orbit equivalence
be a homeomorphism between subsets of full measure, although their
method in fact does seem to yield a $G_{\delta}$.
It is worth  mentioning that in the case of transformations defined on
infinite product spaces  `weakly almost continuous' is equivalent to
`finitary' in the sense of the celebrated work of M.~Keane and
M.~Smorodinsky \cite{KS}.
In a series of subsequent  works this was extended to all odometers \cite{R1},
irrational rotations \cite{R2} and Bratteli-Vershik adic
transformation \cite{HKR} and irreducible Markov chains \cite{RRu}.  All of these
transformations are probability preserving homeomorphisms of compact
metric  spaces.

Recently A.~del Junco and A.~Sahin \cite{dJS}
established an almost continuous version of Dye's theorem:
any two ergodic probability preserving homeomorphisms of  Polish
spaces are almost continuously orbit equivalent. The same is true for
any ergodic homeomorphisms preserving infinite $\sigma$-finite local
measures \cite{dJS}. The purpose of the present paper is to prove an
almost continuous version of Krieger's theorem for systems of type
$III_\lambda$, $0<\lambda\le 1$.

If we only require an orbit equivalence to be a homeomorphism between
invariant sets of full measure then we obtain a weaker notion of
orbit equivalence which we have referred to above here as weakly continuous orbit
equivalence. In fact the main result in \cite{JRW} shows that, under mild assumptions, this is
no stronger than measure-theoretic orbit equivalence. Although the result there
is stated and proved only in the measure-preserving case the statement
and its proof remain valid in the non-singular setting. In fact one
has the following quite general result.

\proclaim{Theorem} Suppose $X$ and $X'$ are separable metric spaces
with Borel probability measures $\mu$ and $\mu'$  and $G,G'$ are
countable ergodic groups of measurable non-singular transformations on $X$ and
$X'$. If
$\phi: X \to X'$ is any measurable non-singular orbit equivalence of
$G,G'$ then there exist elements $\sigma, \tau$ in the measurable
full groups of $G,G'$ such that $\tau \phi \sigma$ (which is also an
orbit equivalence!) is a homeomorphism between invariant subsets
$X'_{0}, Y'_{0}$ of full measure.
\endproclaim

 Section~1 contains the  main definitions. We introduce there a
topological analogue
 $r_{\text{top}}$ of the (measure theoretical) asymptotic range for
the Radon-Nikodym cocycle  of a system.
 It is a closed subgroup of $\Bbb R$ which contains the  asymptotic
range and it is invariant under
 almost continuous orbit equivalence. We give an example of type
$III$ homeomorphisms which are measure-theoretically
orbit equivalent but not almost continuously
orbit equivalent. In Sections~2 and 3  we prove the the main result
of the paper:

\proclaim{Theorem 0.1} Let $(X,\tau,\mu,T)$ and $(X',\tau',\mu', T')$
be
 ergodic non-singular homeomorphisms of Polish spaces.
 If the two systems are either
\roster
\item"\rom{(i)}"
of type $III_\lambda$, $0<\lambda< 1$, and
$r_{\text{top}}(\rho_\mu)=r_{\text{top}}(\rho_{\mu'})=\log\lambda\cdot\Bbb
Z$, or
\item"\rom{(ii)}"
of type $III_1$
\endroster
then they are almost continuously orbit equivalent.
\endproclaim

In the proof we combine techniques  developed in
\cite{dJS} for measure-preserving continuous systems with an approach
to  Krieger's theorem from \cite{KaW} (see also \cite{HO}).
We also achieve a significant technical simplification of the method
of \cite{dJS} by mapping the two spaces involved into a common
symbolic space, thereby removing the need for the many-to-one set
mappings in \cite{dJS}. We do not
know whether the assumption
$r_{\text{top}}(\rho_\mu)= r_{\text{top}}(\rho_\mu')=\log\lambda\cdot\Bbb
Z$ in (i) can be weakened to
$r_{\text{top}}(\rho_\mu)= r_{\text{top}}(\rho_\mu')$.

\head 1. Non-singular homeomorphisms and groups of homeomorphisms of Polish spaces
\endhead

By a dynamical system we mean a quadruple $(X,\tau, \mu,T)$, where
$(X,\tau)$ is a  Polish space, $\mu$ a non-atomic Borel probability
measure
on $X$ and $T$ a $\mu$-non-singular homeomorphism of $X$ satisfying
\roster
\item"$(\bullet)$"
$\mu$ is of full support, i.e. $\mu(O)>0$ for all $O\in\tau$,
\item"$(\bullet)$"
$T$ is ergodic and
\item"$(\bullet)$"
the Radon-Nikodym derivative  $X\ni x\mapsto({d\mu\circ
T}/{d\mu})(x)\in\Bbb R$ is  continuous (i.e. has a continuous version).
\endroster
It follows from these assumptions
that any set of full measure is dense in $X$,  that $X$ is perfect
and that $T$ is topologicaly transitive,
i.e. there exists a point of $X$ whose orbit is dense. The continuous
Radon-Nikodym derivative is defined everywhere in a unique way.
We denote the $T$-orbit equivalence relation by $\Cal R_T$. It is a
Borel
subset of $X\times X$. Let $\rho_\mu:\Cal R_T\to\Bbb R$ stand for the
(logarithm of the) Radon-Nikodym cocycle of $T$, i.e.
$\rho_\mu(x,T^nx)=\log\frac{d\mu\circ T^n}{d\mu}(x)$.

\definition{Definition 1.1}
An invariant dense $G_\delta$-subset of full measure (and the
restriction of
the system to this subset) is called an {\it inessential reduction}
of the
system.

Two systems $(X,\tau, \mu,T)$ and $(Y,\sigma, \nu,S)$ are {\it almost
continuously orbit equivalent} if there are inessential reductions
$X_0\subset X$ and $Y_0\subset Y$ and a homeomorphism $\phi:X_0\to
Y_0$ such that
\roster
\item"(a)"
$\phi(\{T^{n}x: n \in \Bbb Z\}) = \{S^{n}\phi x : n \in \Bbb Z\}$ for all
$x \in X_{0}$,
\item"(b)" $\mu\circ\phi^{-1}\sim\nu$ and the Radon-Nikodym derivative $$
Y_0\ni
y\mapsto{(d\mu\circ\phi^{-1}}/{d\nu)}(y)\in\Bbb R
$$
is (i.e. can be
chosen)
continuous,
\item"(c)"
letting $S' =
\phi^{-1}S \phi$ we have $Tx= {S'}^{n(x)}x$ and $S' = T^{m(x)}x$ where
$n$ and $m$ are continuous on $X_{0}$.
\endroster
\enddefinition

Notice that an inessential reduction of a dynamical system is a
dynamical
system, i.e. the three already mentioned properties marked with
$(\bullet)$ are
satisfied. Since in this paper we study dynamical systems up to
almost
continuous orbit equivalence we will often ignore the difference between a
dynamical system and an inessential reduction of it and thus use the
same
notation for them.
In particular, as
was shown in \cite{dJS}, by
passing to an invariant dense $G_\delta$ subset of full measure we
can and shall assume
without loss of generality that

\roster
\item"$(\bullet)$"
$(X,\tau)$ is 0-dimensional, i.e. there exists a countable base of $X$
consisting of clopen subsets.
\endroster

In the proof of Theorem 0.1 we will establish that the two systems
are almost continuously orbit equivalent to a common intermediate
symbolic system where the action is not by a single homeomorphism but
by a group of homeomorphims so we need to make a few remarks about
almost continuous orbit equivalence in this more general context.
Suppose $G$ and $H$ are countable groups of non-singular homeomorphisms of
Polish probability spaces $X$ and $Y$.
  Note that we
do not assume that
$G$ or $H$ acts freely.
Definition 1.1 of almost continuous orbit equivalence
extends in an obvious way to this context. Condition (c) in that
defintion is replaced  by the following:
letting $G' = \phi^{-1}H\phi$,
for each $h \in
G'$  we have $hx = a(x)x$ where $a: X_{0} \to G$ is continous on
$X_{0}$ and  for each $g \in G$ we have $gx = b(x)x$ where
$b:X_{0} \to G'$ is continuous on $X_{0}$. (Since we have not assumed
freeness the functions $a$ and $b$ here need not be uniquely defined.)
It is clear that the composition of almost continuous orbit
equivalences is again an almost continuous orbit equivalence, after
restriction to a suitable inessential reduction.

If $X$ is a Polish probability space and $P(x)$ is a property of a point
$x \in X$ then we will say that $P(x)$ holds {\it virtually everywhere} if
the set where it holds
contains a dense $G_{\delta}$ subset of full measure.
Any such set will be called virtually full. A function on $X$ will be
called virtually continuous if there is a virtually full subset where
the restriction of the function is continuous. A subset $E$ is virtually
open (closed) if it agrees virtually everywhere with some open (closed) set
and virtually clopen if it is virtually closed and open.
Equivalently,
$E$ is virtually clopen if there is a virtually full subset $X'$ such
that $E \cap X'$ is relatively clopen in $X'$, or $1_{E} $ is
virtually continuous. Since every virtually full set contains an
inessential reduction and we can freely restrict to any inessential
reduction we will often blur the distinction between clopen sets and
virtually clopen sets, as well as between continuous and
virtually continuous functions.

It is also useful to
view almost continuous orbit equivalence in terms of a suitable full
group.
Suppose $G$ is a countable group of non-singular homeomorphisms of a
Polish probability space $X$. We consider all
homeomorphisms $h : X_{0} \to X_{1}$ between virtually full subsets
$X_{0}$ and $X_{1}$ having the form
$h(x) = a(x)x$, where $a:X_{0} \to G$ is continuous on $X_{0}$.
Note that any such $h$ is
necessarily non-singular, with a continuous Jacobian.
We consider two such homeomorphisms to be equivalent if they agree
virtually everywhere.

The collection of equivalence classes is then a
group which we will denote by $[G]_{\text{top}}$ and refer to as the
topological full group of $G$.
Evidently if $\phi$ is an almost continuous orbit equivalence of
groups
$G_{1}$ and $G_{2}$ then $\phi^{-1}[G_{2}]_{\text{top}}\phi =[G_{1}]_{
\text{top}}$.
Conversely it is easy to check that
if $\phi: X_{0} \to Y_{0}$ is any non-singular homeomorphism
between virtually full subsets (invariant or not) which
conjugates the topological full groups then (after restriction to a
suitable inessential
reduction) $\phi$ is an almost continuous orbit equivalence.
(The proof of this is formally the same as the proof of the
corresponding fact in measure-theoretic orbit equivalence.)
Finally we remark that the
groups $G_{i}$ need not consist of everywhere defined homeomorphisms:
it suffices to have virtually everywhere defined  homeomorphisms
whose equivalence classes form a group.

We will also have occasion to use the topological full groupoid
$[[G]]_{\text{top}}$ of
$G$.  We define $[[G]]_{\text{top}}$ as the set of all
(equivalence classes) of homeomorphisms
$h:A \to B$ such that $A$ and $B$ are virtually clopen subsets of $X$ and
$h(x) = a(x)x$ for some $a : A \to G$ which is continuous on $A$.

The following lemma from \cite{dJS} will be used often in this work.

\proclaim{Lemma 1.2}
 Let $A\subset X$ be a clopen subset and
$\mu(A)=\sum_{i=1}^\infty a_i$ with $a_i>0$. Then there exists a
partition
$A=\bigsqcup_{i=1}^\infty A_i$ of $A$ into clopen subsets $A_i$ (on an
inessential reduction) such that $\mu(A_i)=a_i$ for all $i$.
\endproclaim

\demo{Proof}
This is essentially Lemma 2(b) of \cite{dJS} where the $A_{i}$ are only
open. We may take them to be clopen here because of the implicit inessential
reduction. In \cite{dJS} Lemma 2(b) follows from Lemma 2(a) which
asserts that for any $\epsilon > 0$ $A$ may be partitioned into finitely
many clopen sets of measure less than $\epsilon$ and was proved using the
 (measure-preserving) dynamics. However it is easy to give a
 non-dynamical proof. Let $\{U_{i}\}$ be a countable basis for the
 topology and  every $U_i$ is clopen. Then for each $x \in A$ there is a $U_{i}$ of measure less
 than $\epsilon$ which contains $x$. Disjointifying the countably many
 sets so obtained we get a countable partition of $A$ into clopen
 sets of measure less than $\epsilon$ and the existence of a finite
 partition follows immediately. \qed
\enddemo

We recall the concept of the asymptotic range $r(\rho_\mu)\subset\Bbb R$
of
$\rho_\mu$ (\cite{Kr2}, \cite{HO}) and define a topological
counterpart
$r_{\text{top}}(\rho_\mu)$ of this concept.

\definition{Definition 1.3}
\roster
\item"(i)"
A real $t$ belongs to $r(\rho_\mu)$ if for any $\epsilon>0$ and any
subset
$A\subset X$ of positive measure there are $k\ne 0$ and a subset
$B\subset
A$ of positive measure such that $|\rho_\mu(x,T^kx)-t|<\epsilon$ for
all
$x\in B$.
\item"(ii)"
If the above property holds only for all virtually clopen subsets $A$
of positive measure we
write $t\in r_{\text{top}}(\rho_\mu)$. Notice that then the
corresponding
subset $B$ may be chosen virtually clopen.
\endroster
\enddefinition

It is easy to verify that $r_{\text{top}}(\rho_\mu)$ is a closed
subgroup
of $\Bbb R$ and evidently $r_{\text{top}}(\rho_\mu)\supset r(\rho_\mu)$.
 It is also easy to verify that the
asymptotic range and the topological asymptotic range of the
Radon-Nikodym
cocycle are invariant under  measure-theoretic orbit equivalence
and
almost continuous orbit equivalence respectively.

Given a subset $A$ of positive measure, $T$ induces a non-singular
transformation $T_A$ on the space $(A,\mu\restriction A)$. If $A$ is
clopen then $T_{A} \in   [[T]]_{\text{top}}$ and
there exists an inessential reduction $X'$ of $X$ such that
$T_{A \cap X'}$ is a
homeomorphism of $A\cap X'$.  In what follows by $T_A$ we mean $T_{A\cap
X'}$.  Of course, the $T_A$-orbit equivalence relation coincides with
$\Cal R_T\cap (X'\cap A \times X' \cap A)$.

\proclaim{Lemma 1.4} Let $A$ and $B$ be two
non-empty clopen subsets of $X$.
 If one of the following is satisfied:
\roster
\item"\rom{(i)}"
$T$ preserves $\mu$ and $\mu(A)=\mu(B)$ or
\item"\rom{(ii)}"
$T$ is of type $III$
\endroster
then there exist open subsets $A'\subset
A$ and $B'\subset B$ with $\mu(A\setminus A')=\mu(B\setminus B')=0$
and a
homeomorphism $S:A'\to B'$ such that $S \in [[T]]_{\text{ top}}$.
\endproclaim

\demo{Proof} The case (i) was demonstrated  in \cite{dJS}.

We now consider the case (ii). We construct $S$ via an inductive
process. Since $T$ is of
type $III$, there are $N>0$ and mutually disjoint Borel subsets
$(A_i)_{|i|<N}$ of $A$ such that
$0.5\mu(A)<\sum_{|i|<N}\mu(A_i)<\mu(A)$,
the sets $B_i:=T^iA_i$, $|i|<N$, are mutually disjoint Borel subsets
of $B$
and $0.5\mu(B)<\sum_{|i|<N}\mu(B_i)<\mu(B)$ \cite{Kr2}, \cite{HO}.
Fix a small $\epsilon>0$ and
find mutually disjoint clopen subsets $A_i^0\subset A$ such that
$\mu(A_i\triangle A_i^0)<\epsilon\mu(A_i)$, the subsets
$B_i^0:=T^iA_i^0$,
$|i|<N$, are mutually disjoint (clopen) subsets of $B$ and
$\mu(B_i\triangle B_i^0)<\epsilon\mu(B_i)$. Now we can define a
homeomorphism---a `piece' of $S$---from the clopen subset
$A^0:=\bigsqcup_{|i|<N}A_i^0$ of $A$ onto the clopen subset
$B^0:=\bigsqcup_{|i|<N}B_i^0$ of $B$ by setting $Sx=T^ix$ if $x\in
A_i^0$. Notice that $\mu(A^0)>0.5\mu(A)$ and $\mu(B^0)>0.5\mu(B)$
if $\epsilon$ is sufficiently small.
This
completes the first step of the inductive procedure. Now repeat the
same
for the  pair of clopen subsets $A\setminus A^0$ and $B\setminus B^0$
and so on. It remains to concatenate the `pieces' of $S$ that were
defined in all steps to obtain a homeomorphism $\phi \in [[T]]_{\text{top}}$ from an open
subset
$A'$ of $A$, $\mu(A')=\mu(A)$, onto an open subset $B'$ of $B$,
$\mu(B')=\mu(B)$. The graph of $S$ is a subset of $\Cal R_T$ because
each `piece' of $S$ enjoys this property. \qed
\enddemo

\proclaim{Corollary 1.5} If $T$ is of type $III$ and $A$ is
virtually clopen
then
$T_A$ is almost continuously orbit equivalent to $T$.
\endproclaim

For some non-clopen subsets $A\subset X$, the induced system $T_A$
can be well defined too. However, it may not be almost continuously
orbit equivalent to $T$. We illustrate this by the following example.

\example{Example 1.6} Let $C_n:=\{0,1\}$ if $n$ is odd and
$C_n:=\{0,1,\dots,n^2+1\}$ if $n$ is even. Let $T$ stand for the
odometer
on the infinite product space $X=\prod_{n=1}^\infty C_n$ endowed with
the infinite product topology ($C_n$ are all endowed with the discrete
topology). Then $T$ is a homeomorphism of a 0-dimensional compact
metric
space (Cantor set). Fix two reals $0<\lambda,\alpha<1$. We define a
probability measure $\mu_n$ on $C_n$ by setting
$\mu_n(0)=(1+\lambda)^{-1}$
and $\mu_n(1)=\lambda/(1+\lambda)$ for odd $n$ and for even $n$
$$
\mu_n(i)=\cases 1/(n^2(1+\alpha)), &\text{if }i=0\\
\alpha/(n^2(1+\alpha)), &\text{if }i=1\\
n^{-2}, &\text{if }1<i\le n^2 +1.
\endcases
$$
  It is easy to see that the product measure
$\mu=\bigotimes_{n=1}^\infty\mu_n$ on $X$ is non-atomic. Moreover,
$T$ is
$\mu$-non-singular and ergodic. Let $A:=C_1\times C_2'\times C_3\times
C_4'\times\dots$, where $C_n':=\{2,3,\dots, n^2+1\}$. Then $A$ is a
compact
meager subset of $X$ and $\mu(A)>0$. It is easy to see that $T_A$ is
a well defined
homeomorphism of $A$, in fact an odometer, and the dynamical
system
$(T_A, \mu\restriction A)$ is of type $III_\lambda$. Hence $T$ is
also of
type $III_\lambda$. Notice that the map $x\mapsto\rho_\mu(x,Tx)$ is
continuous everywhere on $X$ except of a singleton. It is easy to
verify
that $\{\log\lambda, \log\alpha\}\subset r_{\text{top}}(\rho_\mu)$.
Since
$\rho_\mu$ takes values in the group generated by $\log\lambda$ and
$\log\alpha$, it follows that $r_{\text{top}}(\rho_\mu)$ is the
closure of
this group. Thus $r_{\text{top}}(\rho_\mu)\ne r(\rho_\mu)$ whenever
$\alpha$ is not a power of $\lambda$. It follows that $T$ and $T_A$
while being orbit
equivalent are not almost continuously orbit equivalent.
\endexample

\proclaim{Lemma 1.7} If $r_{\text{top}}(\rho_\mu)\cap[a,b]=\emptyset$
then there is a clopen subset
$A\subset X$ with $\rho_\mu(x,y)\not\in[a,b]$ for all $(x,y)\in\Cal
R_T\cap(A\times A)$.
\endproclaim
\demo{Proof} We first find for each $t\in[a,b]$, a clopen subset
$B_t\subset X$ and $\epsilon_t>0$ such that
$\rho_\mu(x,y)\notin(t-\epsilon_t,t+\epsilon_t)$ for all $\Cal
R_T\cap(B_t\times B_t)$. Select a finite collection of the intervals
$(t-\epsilon_t/2,t+\epsilon_t/2)$ that cover the entire $[a,b]$ and
enumerate them as $(t_i-\epsilon_i/2,t_i+\epsilon_i/2)$, $1\le i\le
l$,
with corresponding sets $B_i$. Since $T$ is ergodic and the
Radon-Nikodym
derivative $d\mu\circ T/{d\mu}$ is continuous, there exist a clopen
subset
$C_2\subset B_1$ and $k\in\Bbb Z$ with $T^kC_2\subset B_2$ and the map
$C_2\ni x\mapsto\rho_\mu(x,T^kx)$ is close to a constant. Take
$(x,y)\in\Cal R_T\cap(C_2\times C_2)$. Then
$$
\align
\rho_\mu(x,y) &\notin(t_1-\epsilon_1,t_1+\epsilon_1)\text{ and}\\
\rho_\mu(x,y)
&=\rho_\mu(x,T^kx)+\rho_\mu(T^kx,T^ky)-\rho_\mu(y,T^ky)\notin
(t_2-\epsilon_2,t_2+\epsilon_2).
\endalign
$$
Continue in this fashion to get clopen sets $C_2\supset
C_3\supset\cdots\supset C_l$ such that for some $k_i$ we have
$T^{k_i}C_i\subset B_i$ and the map $C_i\ni
x\mapsto\rho_\mu(x,T^{k_i}x)$
is close to a constant. Setting $A:=C_l$ proves the lemma. \qed
\enddemo

\proclaim{Lemma 1.8} There is an inessential reduction
$\widetilde X\subset X$ and a countable base
$(O_i)_{i=1}^\infty$ of the topology in $\Bbb R$ such that the sets
$$
\{x\in\widetilde X\mid\rho_\mu(x,Tx)\in O_i\}, \qquad i=1,2, \ldots
$$
are clopen in
$\widetilde X$
(in the relative topology).
\endproclaim
\demo{Proof} We first notice that there are no more than countably
many
$t\in\Bbb R$, say bad $t$, such that the closed subsets
$$
F_t:=\{x\in\widetilde X\mid\rho(x,Tx)=t\}
$$
have positive measure.
Chose a
countable family of intervals $(O_i)_{i=1}^\infty$ in $\Bbb R$ which
generate the topology and whose end points are not bad. Let $F$ be the
smallest $T$-invariant subset that contains $F_t$ when $t$ runs
through the
end
points of all $O_i$, $i=1,2,\dots$. It remains to put $\widetilde
X:=X\setminus F$. \qed
\enddemo

\definition{Definition 1.9} Given $r>0$, we let
$I_r:=\{0,1,\dots,r-1\}$ and
endow it with the discrete topology. A subrelation $\Cal S\subset\Cal
R_T$
is called $r$-{\it uniform} if there are a clopen subset $B\subset X$
 and a homeomorphism $h:I_r\times B\to X$ such that
$h(I_r\times\{b\})=\Cal S(b)$ and $h(0,b)=b$ for $b\in B$ and for each
$i$ the map
$$
B \ni h(0,b) \mapsto h(i,b) \in h(\{i\} \times B)
$$
belongs to
$[[T]]_{   \text{top}}$.
  We call $B$ an $\Cal
S$-{\it fundamental subset} and $h$ an $(\Cal S,B)$-{\it splitting
homeomorphism}. Associated to $\Cal S$ we have a finite subgroup
$G_{\Cal S}$ of $[T]_{\text{top}}$ corresponding to the action of the
permutation group $S_{r}$ on $I_{r}$. We also have a natural
projection $p: X \to B$ mapping $h(i,b)$ to $h(0,b)$ for each $i$.
\enddefinition

It is worth noting that neither $B$ nor $h$ is determined
uniquely by
$\Cal S$ or $(\Cal S,B)$ respectively. In a similar way, $G_{\Cal S}$ depends really on $(\Cal S, B,h)$ and not only on $\Cal S$ alone. Each $q$-uniform subrelation
$\Cal
T$ of $\Cal R_T\cap (B\times B)$ extends canonically, via $\Cal S$, to an
$rq$-uniform
subrelation
$$
\widetilde{\Cal T}:=\{(x,y)\in \Cal R_T\mid (p(x),p(y))\in\Cal T\}
$$
of $\Cal R_T$. We call it the {\it natural extension} of $\Cal T$ to
$X$.
Each $\Cal T$-fundamental subset  is also a $\widetilde{\Cal
T}$-fundamental subset.

The following topological version of the non-singular Rokhlin lemma
can be
proved in the same way as its measure-preserving counterpart (Lemma
6, \cite{dJS}).

\proclaim{Lemma 1.10} Given $\epsilon>0$ and an integer $N>0$, there is a
clopen subset $A\subset X$ such that the subsets $T^iA$, $i=0,\dots,
N$,
are mutually disjoint and
$
\mu(X\setminus\bigsqcup_{i=0}^NT^iA)<\epsilon.
$
\endproclaim

We need this lemma to obtain the following result.

\proclaim{Lemma 1.11} Let $\Cal S$ be a uniform subrelation of $\Cal
R_T$
and let $B\subset X$ be an $\Cal S$-fundamental subset. For each
$\epsilon>0$, there exists a uniform subrelation $\Cal T\subset\Cal
R_T\cap
(B\times B)$ such that $\mu(\{x\in X\mid Tx\in \widetilde{\Cal
T}(x)\})>1-2\epsilon$.
\endproclaim

\demo{Proof}
Take $\delta>0$ such that whenever $C$ is a subset of $B$ with
$\mu(C)<\delta\mu(B)$ then $\mu(\bigcup_{x\in C}\Cal S(x))<\epsilon$.

Denote by $T_B:B\to B$ the induced homeomorphism. Select a
clopen subset $A\subset B$ and a positive integer $M$ such that
$\mu(A)>(1-\delta/2)\mu(B)$ and $T \Cal S(x)
\subset\bigsqcup_{|i|<M}\Cal
S(T_B^ix)$ for each $x\in A$. We deduce from the Lemma~1.10 that for
a sufficiently large
$N>M$, there exists a clopen subset $B_0\subset B$ such that the
subsets
$T_B^iB_0$, $i=0,\dots,N-1$ are mutually disjoint and
$$
\mu\bigg(B\setminus\bigsqcup_{i=M}^{N-M}T_B^iB_0\bigg)<\delta\mu(B)/2.
$$
Partition
the subset $B\setminus\bigsqcup_{i=0}^{N-1}T_B^i B_0$ into $N$
non-empty clopen
subsets $D_j$, $0\le j<N$. By Lemma~1.4(ii), there are homeomorphisms
$\gamma_j:D_0\to D_j$ belonging to $[[T]]_{\text{top}}$. Now we
define an equivalence relation ${\Cal T}$ on $B$ by setting
$(x,y)\in{\Cal
T}$ if either
\roster
\item"---"
$x=T_B^iz$ and $y=T_B^jz$ for some $z\in B_0$ and $0\le i,j<N$ or
\item"---"
$x=\gamma_iz$ and $y=\gamma_jz$ for some $z\in D_0$ and $0\le i,j<N$.
\endroster
It is obvious that ${\Cal T}$ is an $N$-uniform subrelation of $\Cal
R_T\cap(B\times B)$. Let $E$ stand for the union of all subsets $\Cal
S(x)$
when $x$ runs through $A\cap\bigsqcup_{M<i<N-M}T^i_BB_0$. Then
$\mu(E)>1-\epsilon$
and $Tx\in\widetilde{\Cal T}(x)$ for all $x\in E$. \qed
\enddemo

We remark that the set
$O := \{x\in X\mid Tx\in \widetilde{\Cal
T}(x)\}$ in the conclusion of Lemma~1.11 is automatically (virtually) clopen
 and that $T|_{O} \in [[G_{ \widetilde{\Cal
T}(x)}]]_{\text{top}}$.

\head 2. Type $III_\lambda$ with $0<\lambda<1$
\endhead

In this section we prove Theorem~0.1 in the case when $T$ is of
type $III_\lambda$ with $0<\lambda<1$.

\proclaim{Proposition 2.1} Let $(X,\tau,\mu,T)$ be a dynamical
system  with $r_{\text{top}}(\rho_\mu)=\log\lambda\cdot\Bbb Z$.
(We do not assume that $r(\rho_\mu)=\log\lambda\cdot\Bbb Z$).
Then
there exists an equivalent probability measure $\widetilde\mu$ such
that
$\rho_{\widetilde\mu}(x,y)\in\log\lambda\cdot\Bbb Z$ for all
$(x,y)\in\Cal
R_T$ and $d\widetilde\mu/d\mu$ is continuous (when restricted to an
inessential reduction).
\endproclaim
\demo{Proof}
We first note that if there are $\eta>0$ and a clopen
subset
$A\subset X$ such that
$$
\rho_\mu(x,y)\not\in[\log\lambda+\eta,-\eta]
\text{ for
all }(x,y)\in\Cal R_T\cap (A\times A)
$$
then
$\rho_\mu(x,y)\in\log\lambda\cdot\Bbb Z+(-\eta,\eta)$ for all
$(x,y)\in\Cal
R_T\cap (A\times A)$.

Fix a summable sequence $(\eta_{i})_{i=0}^\infty$ of positive reals.
It follows from
Lemma~1.7 and the last remark that there exists a  clopen subset
$A_0 \subset X $ such that
$$
\rho_\mu(x,y)\in\log\lambda\cdot\Bbb Z+(-\eta_1,\eta_1)
\text{
for all }
(x,y)\in\Cal R_T\cap (A_0 \times A_0).
$$
Passing
to an inessential reduction, we can represent $X$ as
$X=\bigsqcup_{k=0}^\infty A_k$, where $A_k$ is clopen and
$T^{-k}A_k \subset A_0$. We now apply Lemma~1.8 to find a
continuous function $\phi_1:X\to(\log\lambda,-\log\lambda)$ such that
$\phi_1=0$ on $A_0^{(i)}$ and
$$
\phi_1(x)+\rho_\mu(x,T^{-k}x)\in\log\lambda\cdot\Bbb
Z+(-\eta_1,\eta_1)
\text{ for all }x\in A_k.
$$
 Let $\mu_1:=\exp(\phi_{1})\cdot\mu$. Then it
is
easy to verify that
$$
\rho_{\mu_1}(x,y)\in \log\lambda\cdot\Bbb Z+(-3\eta_1,3\eta_1)\text{
for
all }(x,y)\in\Cal R_T.   \tag{2-1}
$$
We now repeat the process with $\mu_{1}$ in place of $\mu$ to obtain
$\mu_{2} = \exp(\phi_{2})\mu_{1}$ such that
$$
\rho_{\mu_2}(x,y)\in \log\lambda\cdot\Bbb Z+(-3\eta_2,3\eta_2)\text{
for
all }(x,y)\in\Cal R_T.
$$
However, now because of \thetag{2-1} we see that $\phi_{2}$ may be chosen
with values in $(-3\eta_1,3\eta_1)$.
Continuing in this fashion we obtain a sequence $(\mu_{i})_{i=1}^\infty$ of measures
such that $\mu_{i+1} = \exp(\phi_{i})\mu_{i}$,
$\phi_{i}$ is continuous,   $|\phi_{i+1}| < 3\eta_{i}$
and
$$
\rho_{\mu_i}(x,y)\in \log\lambda\cdot\Bbb Z+(-3\eta_i,3\eta_i)\text{
for
all }(x,y)\in\Cal R_T.
$$
The series $\sum_{i>0} \phi_{i}$ converges uniformly to a continuous
function $\phi$ and we have only to set $\widetilde \mu = \exp(\phi)\mu$ to
complete the proof.
We note that  $\widetilde\mu$ (as well as $\widetilde\mu_i$) can be
infinite  by construction. Nevertheless it is locally finite. Hence,
after multiplying by a suitable continuous function taking values in
$\{n\log \lambda  :n \in \Bbb Z\}$,
  we may assume without loss of
generality that $\widetilde\mu$ is finite.
\qed
\enddemo

We will call $\widetilde\mu$ a {\it special measure} for $T$.

\remark{Remark \rom{2.2}}If
$r_{\text{top}}(\rho_\mu)=\log\lambda\cdot\Bbb Z$ and $\widetilde
\mu$ is special then the ``measure-preserving''
subrelation $\Cal S:=\{(x,y)\in\Cal
R_T\mid
\rho_{\widetilde\mu}(x,y)=0\}$ of $\Cal R_T$ is ergodic
if and only if
$r(\rho_\mu)=\log\lambda\cdot\Bbb Z$.

To see this let us agree to say that measurable subsets $A$ and $B$
communicate via $l\in\Bbb Z$ if there are
subsets $A' \subset A,\ B' \subset B$ of positive measure
and $k \in \Bbb Z$ such that $T^{k}A' = B'$ and
$\rho_{\widetilde \mu}(x, T^{k}x) = l \log \lambda$ for all $x \in A'$.

Now if $r(\rho_{\mu})=\log\lambda\cdot\Bbb Z$ and
 $A$ and $B$ are measurable sets of
positive measure then certainly $A$ and  $B$ communicate via some
$l$.
Since $T$ is type
$III_{\lambda}$ and $\widetilde \mu$ is special we also have that
every subset of $B$  of positive measure communicates with itself via $-l$ and it follows that $A$ and $B$ communicate via $0$, as desired.

On the other hand if $\Cal S$ is ergodic, since
$r_{\text{top}}(\rho_\mu)=\log\lambda\cdot\Bbb Z =
r_{\text{top}}(\rho_{\widetilde{\mu}})$ then there certainly exist
(clopen)
$A$ and $B$ which communicate via $1$. If $C$ is an arbitrary
measurable set of positive measure then by ergodicity of $\Cal S$
each subset of $C$ of positive measure communicates with both $A$ and $B$ via $0$ and
`by composition' it follows that $C$ communicates with $C$ via $1$,
showing that $\lambda \in r(\rho_{\widetilde \mu}) = r(\rho_{ \mu}) $.
Thus $r(\rho_{ \mu}) \supset \log\lambda\cdot\Bbb Z$ and we already have the
opposite inclusion.
\endremark

The following statement is a refinement of Lemma~1.4(ii) for type
$III_\lambda$ systems.

\proclaim{Lemma 2.3} Let $T$ be  of type $III_\lambda$,
$r_{\text{top}}(\rho_\mu)=\log\lambda\cdot\Bbb Z$ and $\mu$ a special
measure. Given two clopen subsets $A,B\subset X$ with
$\mu(B)=\lambda^{k}\mu(A)$, there exists a homeomorphism $S:A\to B$
such
that $S \in [[T]]_{\text{top}}$
  and $\rho_\mu(x,Sx)=k\log \lambda$ for all
$x\in A$.
\endproclaim

\demo{Proof} Let $n(x):=\min\{n>0\mid\rho_\mu(x,T^nx)=0\}$ and
$Rx:=T^{n(x)}x$.
Since $\rho_{\mu}$ is continuous and takes values in the discrete
set $\log\lambda\cdot\Bbb Z$, it is easy to see that $n$ is virtually continuous
so that $R \in [[T]]_{\text{top}}$.
The $R$-orbit equivalence relation co-incides with
 $\Cal S$, the measure-preserving subrelation. Hence  $R$ is ergodic by Remark~2.2. Since $T$ is ergodic of type
$III_\lambda$, there are a
non-empty clopen subset $A_0\subset X$ and $l\in\Bbb Z$ such that
 $\rho_\mu(x,T^lx)=k\log \lambda$ for all $x\in A_0$.
 In view of Lemma~1.2 we may assume that
 $\mu(A_{0})=N^{-1}\mu(A)$ for some $N\in\Bbb Z$. Again by Lemma~1.2,
 there are
 partitions $A=\bigsqcup_{i=1}^N A_i$ and $B=\bigsqcup_{i=1}^N B_i$
of $A$
 and $B$ into clopen subsets $A_i$ and $B_i$ such that
 $\mu(A_i)=N^{-1}\mu(A_0)$ and $\mu(B_i)=N^{-1}\mu(B_0)$ for all $i$.
By
 Lemma~1.4(i), there exist homeomorphisms $\gamma_i:A_i\to A_{0}$ and
 $\delta_i:B_i\to T^{l}A_{0}$ with $\gamma_{i}, \delta_{i} \in
 [[R]]_{\text{top}}$.
It remains to set $Sx:=\delta_i^{-1}T^{l}\gamma_ix$ for all $x\in
A_i$,
$i=1,\dots,N$. \qed
\enddemo

\proclaim{Lemma 2.4} Let $(\mu, T)$ be as in Lemma~2.3. Let
$D,\widetilde
D$ be two countable discrete spaces and $B$ a clopen subset of $X$.
Let
$\nu$ be a finite measure on $I_r\times \widetilde D$ such that
$\nu(i,d)/\nu(j,d)$ is a power of $\lambda$ for all $i,j\in I_r$ and
$d\in\widetilde D$. Given continuous onto maps $u:B\to D$ and
$\widetilde
v:I_r\times\widetilde D\to D$ such that $\mu\circ
u^{-1}=\nu\circ\widetilde v^{-1}$, there are an $r$-uniform
equivalence relation $\Cal S\subset\Cal R_T$ on $B$, an $\Cal
S$-fundamental subset $\widetilde B\subset B$, an $(\Cal S,\widetilde
B)$-splitting homeomorphism $v:I_r\times\widetilde B\to B$ and a
measure-preserving continuous onto map $\widetilde u:\widetilde
B\to\widetilde D$
such that the following diagram commutes
$$
\CD
 (B,\mu) @<{v}<<    (I_r\times\widetilde B,\mu\circ v)\\
 @V{u}VV    @VV{\text{id}\times\widetilde u}V\\
 (D,\nu\circ \widetilde v^{-1}) @<{\widetilde v}<<
(I_r\times\widetilde D,
\nu )
\endCD
$$
and $\rho_\mu(v(j,\widetilde b),v(i,\widetilde
b))=\log(\nu(j,\widetilde u(\widetilde b))/\nu(i,\widetilde
u(\widetilde b)))$ for all $i,j\in I_r$ and $\widetilde
b\in\widetilde B$.
Here  $\widetilde B$ and $\widetilde D$ are identified with the
`bottom' levels $\{0\}\times \widetilde B$ and $\{0\}\times
\widetilde D$ of $I_r\times \widetilde B$ and $I_r\times \widetilde
D$ respectively equipped with the induced measures.
\endproclaim

\demo{Proof} Since $ \mu(u^{-1}(d))=\nu(\widetilde v^{-1}(d)) \text{
for each } d\in D, $ it follows from Lem\-ma~1.2 that there exists a
partition $B=\bigsqcup_{(i,\widetilde d)\in I_r\times\widetilde
D}B_{i,\widetilde d}$ of $B$ into clopen subsets such that
$\mu(B_{i,\widetilde d})=\nu(i,\widetilde d)$ and $u(B_{i,\widetilde
d})=\widetilde v(i,\widetilde d)$ for all $(i,\widetilde d)\in
I_r\times\widetilde D$. Hence the ratio $\mu(B_{i,\widetilde
d})/\mu(B_{0,\widetilde
d})$ is a power of $\lambda$ and then Lemma~2.3 yields that there
exists a
homeomorphism $\gamma_{i,\widetilde d}:B_{i,\widetilde d}\to
B_{0,\widetilde d}$  such that $\gamma_{i,\widetilde d}\in[[T]]_{\text{top}}$ and 
$$
\rho_\mu(x,\gamma_{i,\widetilde
d}x)=\log(\mu(B_{0,\widetilde d})/\mu(B_{i,\widetilde
d}))\text{ for all }x\in B_{i,\widetilde d}\,.
 $$
Let $\Cal S$ be the equivalence
relation
generated by the graphs of all $\gamma_{i,\widetilde d}$, $i\in I_r$,
$\widetilde d\in\widetilde D$. It is $r$-uniform and the subset
$\widetilde
B:=\bigsqcup_{d\in\widetilde D}B_{0,\widetilde d}\subset B$ is $\Cal
S$-fundamental. It remains to put $\widetilde u=\widetilde d$ on
$B_{0,\widetilde d}$
and $v(i,b)=\gamma_{i,\widetilde d}^{-1}b$ if $(i,b)\in I_r\times
B_{0,\widetilde d}$. \qed
\enddemo

Now we state and prove the main result of this section---part (i) of
Theorem~0.1.

\proclaim{Theorem 2.5} Let $(X,\tau,\mu,T)$ and $(X',\tau',\mu', T')$
be
two  dynamical systems of type $III_\lambda$, $0<\lambda<1$, and
$r_{\text{top}}(\rho_\mu)=r_{\text{top}}(\rho_{\mu'})=\log\lambda\cdot\Bbb
Z$. Then these systems are almost continuously orbit equivalent.
\endproclaim
\demo{Proof} By Proposition~2.1 we may assume without loss of
generality
that $\mu$ and $\mu'$ are special for $T$ and $T'$ respectively. Let
$m$
and $m'$ stand for the complete metrics compatible with the
topologies on
$X$ and $X'$ respectively. Fix a sequence $\epsilon_n\to 0$. Suppose
that
we have inessential reductions $X_0\subset X$ and $X_0'\subset X'$, a
sequence of positive integers $r_n$ and two nested sequences of clopen
subsets
$$
X_0=:B_0\supset B_1\supset\cdots, \ \ X_0'=:B_0'\supset
B_1'\supset\cdots
$$
in $X_0$ and $X_0'$ respectively such that $B_n$ is a fundamental
subset of
an $r_n$-uniform equivalence relation $\Cal R_n\subset\Cal
R_T\cap(B_{n-1}\times B_{n-1})$ and $B_n'$ is a fundamental subset of
an
$r_n$-uniform equivalence relation $\Cal R_n'\subset\Cal
R_{T'}\cap(B_{n-1}'\times B_{n-1}')$. Suppose also that we can choose
for
any $n$,
\roster
\item"$\bullet$"
an $(\Cal R_n, B_n)$-splitting homeomorphism $h_n:I_{r_n}\times B_n\to
B_{n-1}$,
\item"$\bullet$"
a measure preserving continuous map $p_n$ of $B_n$ onto a discrete
countable set $D_n$,
\item"$\bullet$"
a measure preserving  onto map $g_n:I_{r_n}\times D_n\to D_{n-1}$,
\item"$\bullet$"
similar objects $h_n', p_n', g_n'$ for the second system
\endroster
in such a way that the infinite diagram
$$
\CD
 X_0 @<h_1<<  I_{r_1}\times B_1 @<{\text{id}\times h_2}<<
I_{r_1}\times I_{r_2}\times B_2 @<{\text{id}\times\text{id}\times
h_3}<<
\cdots\\
 {} @. @V{\text{id}\times p_1}VV
 @V{\text{id}\times\text{id}\times p_2}VV\\
 {} @. I_{r_1}\times D_1 @<{\text{id}\times g_2}<<  I_{r_1}\times
I_{r_2}\times D_2 @<{\text{id}\times\text{id}\times g_3}<<\cdots\\
@. @A{\text{id}\times p_1'}AA
 @A{\text{id}\times\text{id}\times p_2'}AA\\
 X_0' @<h_1'<<  I_{r_1}\times B_1' @<{\text{id}\times h_2'}<<
I_{r_1}\times I_{r_2}\times B_2' @<{\text{id}\times\text{id}\times
h_3'}<<
\cdots
\endCD
\tag2-2
$$
commutes.

Since each arrow in the upper line is a homeomorphism the measure
$\mu$ in $X_{0}$ induces a measure on each space in the upper line
and
these measures in turn project onto the spaces in the middle line.
We also require of \thetag{2-2} that these measures on the middle
line co-incide with the measures coming in the same way from the
lower line. Thus every space in the picture carries a measure and the
arrows are all measure-preserving.
 For
any $n$, we have a continuous onto map from $X_0$ to
$I_{r_1}\times\cdots\times
I_{r_n}\times D_n$ which is a composition of $n$ reversed horizontal
arrows
plus a vertical arrow. We denote this map by $\pi_n$ and view it as a
partition of $X_0$
into
clopen subsets. We will denote by
$\widetilde {\Cal R}_n$  the natural
extension of $\Cal R_n$ to $X_0$.

We further require of  \thetag{2-2}  that
\roster
\item"$(\triangle)$"
 for each odd $n$, the diameter of (every atom of) $\pi_n$ is less than $\epsilon_n$ with
respect to $m$;

\item"$(\diamondsuit)$" for each odd $n$ there is a clopen subset $O_n\subset X_0$ such
that
$\mu(O_n)>1-\epsilon_n$ and $Tx\in\widetilde{\Cal R}_n(x)$ for all $x\in
O_n$.

\item"$(\bigtriangledown)$" for all $n$ the Radon-Nikodym derivative on ${\Cal
R}_n$ is equivariant under $p_n$, i.e. given $d\in D_n$ and $i,j\in
I_{r_n}$, the map $p_n^{-1}(d)\ni
b\mapsto\rho_\mu(h_n(i,b),h_n(j,b))$ is constant.

\item"$(\heartsuit)$" similar conditions hold for the `bottom line'
of \thetag{2-1} but with ``odd'' replaced by ``even'' $n$ in $(\triangle)$ and
$(\diamondsuit)$.
\endroster

Under all these assumptions on \thetag{2-2} we proceed to construct our orbit
equivalence.
  Denote by $Z$ the topological and measure-theoretic inverse limit of
the
middle line in \thetag{2-2} and denote by $\nu$ the measure on $Z$.
It is obvious that $Z$ is a
0-dimensional Polish space. An element $z$ of $Z$ is a sequence
$(z_n)_{n=1}^\infty$ such that $z_n=(i_1,\dots,i_n,d_n)$ for some
$i_k\in
I_{r_k}$, $1\le k\le n$, and $d_n\in D_n$ with
$d_n=g_{n+1}(i_{n+1},d_{n+1})$, $n\ge 1$.  Since
\thetag{2-2}
commutes and all the arrows are onto, there is a continuous map
$\phi:X_0\to Z$ whose image is dense in $Z$. It follows from~$(\triangle)$
that $\phi$ is one-to-one and open as a map from $X_0$ to $\phi(X_0)$
endowed with the induced topology, that is $\phi$ is a homeomorphism
onto its image. Since $X_0$ is Polish, $\phi(X_0)$
is
Polish. It follows by the Alexandrov-Hopf theorem that $\phi(X_0)$ is a $G_\delta$ in
$Z$. Note also that $\mu \circ \phi^{-1} = \nu$.

We let $G$ denote the countable subgroup of $[T]_{\text{top}}$
which is the increasing
union of the finite subgroups $G_{n}:= G_{\widetilde{\Cal R}_n}$. We claim that,
because of~$(\diamondsuit)$, $T \in  [G_{}]_{\text{top}}$, so that
in fact $[T]_{\text{top}} = [G_{}]_{\text{top}}$. To see this note
that for $n$ odd we have $T|_{O_{n}} \in [[G_{n}]]_{\text{top}}$, that is
$Tx = g_{n}(x)(x)$ where $g_{n} : O_{n} \to G_{n}$ is virtually
continuous. (See the remark after the proof of Lemma 1.11.)
Let $O = \bigcup_{n \geq 0} O_{2n+1}$, a
virtually open set of full measure. If we define $g: O \to G$ by
$g(x) = g_{2n+1}(x) $ for $x \in O_{2n+1} \backslash O_{2n-1} $
 then $g$ is virtually continuous, since each
$O_{2n+1} \backslash O_{2n-1}$ is clopen and $Tx = g(x)x$ for $x \in O$.
This shows that $T \in  [G_{}]_{\text{top}} $ as claimed.

The permutation group $S_{r_{1}r_{2} \ldots r_{n}}$ acts on
$I_{r_{1}} \times \dots \times
I_{r_{n}¥}¥$ and these permutations lift to homeomorphisms of $Z$,
comprising a group which we denote $H_{n}$. It is clear that $\phi
G_{n}\phi^{-1} = H_{n}$ so $\phi G \phi^{-1} = H$ where $H = \bigcup_{n}
H_{n}$. Thus we have
 $$
\phi[T]_{\text{top}}\phi^{-1} = \phi[G]_{\text{top}}\phi^{-1} =
    [H]_{\text {top}}.
   $$

 In a similar way---using
$(\heartsuit)$---one can define a measure-preserving map
$\phi':X_0'\to Z$ which is a homeomorphism onto a dense $G_{\delta}$ such that
 $$
 \phi'[T]_{\text{top}}(\phi')^{-1} =     [H]_{\text {top}}
$$
and $\mu' \circ (\phi')^{-1} = \nu$.
Now the map $\psi = (\phi')^{-1}\circ\phi$ is virtually everywhere defined and
$\psi[T]_{\text{top}}\psi^{-1} = [S]_{\text{top}} $ so $\psi$ is a
(measure-preserving) almost continuous orbit equivalence.

It remains to explain how to construct \thetag{2-2} satisfying
$(\triangle)$, $(\diamondsuit)$, $(\bigtriangledown)$ and
$(\heartsuit)$. This will be done via
an inductive process. Suppose we have already
constructed~\thetag{2-2}
right up to an odd index $n$. Moreover, suppose that the `bottom
square' of
\thetag{2-2} with index $n+1$ also has been constructed. Thus the
following
fragment of \thetag{2-2} commutes:
$$
\CD
B_n @.\\
@V{p_n}VV @.\\
D_n @<g_{n+1}<< I_{r_{n+1}}\times D_{n+1}\\
@A{p_n'}AA @AA{\text{id}\times p_{n+1}'}A\\
B_n' @<h_{n+1}'<< I_{r_{n+1}}\times B_{n+1}'\\
\endCD
\tag2-3
$$
Recall that the arrows are assumed measure preserving.
Since $\mu'$ is special and $(\heartsuit)$ and $(\bigtriangledown)$
hold,
 it follows that  the measure on $I_{r_{n+1}}\times D_{n+1}$, say
$\nu_{n+1}$, is such that $\nu_{n+1}(i,d)/\nu_{n+1}(j,d)$ is a power
of $\lambda$ for all $i,j\in I_{r_{n+1}}$ and $d\in D_{n+1}$.
Then we can apply Lemma~2.4 to
find  an $r_{n+1}$-uniform subrelation $\Cal
R_{n+1}\subset\Cal R_T\cap(B_n\times B_n)$, an $\Cal
R_{n+1}$-fundamental
subset $B_{n+1}\subset B_n$, a $(\Cal R_{n+1}, B_{n+1})$-splitting
homeomorphism $h_{n+1}:I_{r_{n+1}}\times B_{n+1}\to B_n$ and a measure
preserving continuous onto map $p_{n+1}:B_{n+1}\to D_{n+1}$ such that
the
diagram
$$
\CD B_n @<{h_{n+1}}<< I_{r_{n+1}}\times B_{n+1}\\
@Vp_nVV @VV{\text{id}\times p_{n+1}}V\\
D_n @<{g_{n+1}}<< I_{r_{n+1}}\times D_{n+1}
\endCD
$$
commutes and the Radon-Nikodym derivative on $\Cal
R_{n+1}$ is equivariant under $p_{n+1}$ (i.e. $(\bigtriangledown)$
holds for $n+1$). Now apply Lemma~1.11 to find $r_{n+2}>0$, an
$r_{n+2}$-uniform
subrelation $\Cal R_{n+2}\subset\Cal R_T\cap(B_{n+1}\times B_{n+1})$
and a
subset $O_{n+2}\subset X_0$ satisfying $(\diamondsuit)$. Select an
$\Cal
R_{n+2}$-fundamental subset $B_{n+2}\subset B_{n+1}$ and an $(\Cal
R_{n+2},
B_{n+2})$-splitting homeomorphism $h_{n+2}:I_{r_{n+2}}\times
B_{n+2}\to
B_{n+1}$. Take now any countable onto map $p_{n+2}$ from $B_{n+2}$ to
a
countable discrete set $D_{n+2}$ such that the diameter of the
corresponding partition $\pi_{n+2}$ of $X_0$ is less then
$\epsilon_{n+2}$,
i.e. $(\triangle)$ is satisfied.
`Refining' $p_{n+2}$, if necessary, we can also assume that there
exists an onto map $g_{n+2}:I_{r_{n+2}}\times D_{n+2}\to D_{n+1}$
such that the
 diagram
$$
\CD
B_{n+1} @<h_{n+2}<< I_{r_{n+2}}\times B_{n+2}\\
@V{p_{n+1}}VV @VV{\text{id}\times p_{n+2}}V\\
D_{n+1} @<g_{n+2}<< I_{r_{n+2}}\times D_{n+2}\\
@A{p_{n+1}'}AA @.\\
B_{n+1}' @.
\endCD
\tag2-4
$$
commutes and
$(\bigtriangledown)$ holds for $n+2$.
The diagram~\thetag{2-4} is a `reversed' version of \thetag{2-3} but
with index $n+2$.
Iterating this argument infinitely many times we come to~\thetag{2-2}.
For that $(\bigtriangledown)$ is used: constructing \thetag{2-2}
from `measure preserving' diagrams \thetag{2-3} and \thetag{2-4} we
obtain a fragment of \thetag{2-2} which is also measure preserving.

It should be noted that every time when we apply Lemmata 2.4 and 1.11
to
define $B_{n+1}$ or $B_{n+1}'$ we change the previously defined
subsets
$B_{1},\dots, B_n$ and $B_{1}',\dots, B_n'$. Indeed, we reduce them
to some
inessential reductions of $X$ and $X'$. Making countably many steps we
intersect countably many such inessential reductions to obtain $X_0$
and
$X_0'$ after all steps done. \qed

\enddemo

\head 3. Type $III_1$
\endhead

In this section we prove Theorem~0.1 in the case when $T$ is of
type $III_1$. We start with a crucial lemma which refines
Lemma~1.4(ii) (cf. Lemma~2.3).

\proclaim{Lemma~3.1} Let $T$ be  of type $III_1$ and let $A,B$ be two
non-empty clopen subsets in $X$. Then for each $\epsilon>0$, there are
open subsets $A'\subset A$, $B'\subset B$ and a homeomorphism
$S:A'\to B'$ such that $\mu(A\setminus A')=\mu(B\setminus B')=0$,
$(x, Sx)\in\Cal R_T$ for all $x\in A'$, the map $A'\ni
x\to\rho_\mu(x,Sx)$ is continuous and
$$
|\rho_\mu(x,Sx)-\log(\mu(B)/\mu(A))|<\epsilon\text{ for all }x\in A'.
$$
\endproclaim

\demo{Proof} We construct $S$ via an inductive process. Since the
real $w:=\log(\mu(B)/\mu(A))$ belongs to $r(\rho_\mu)$, it follows
from \cite{Kr1}, \cite{HO} that there are $N>0$ and mutually disjoint
Borel subsets
$(A_i)_{|i|<N}$ of $A$ such that
\roster
\item"$(\circ)$"
$1/2\mu(A)<\sum_{|i|<N}\mu(A_i)<3/4\mu(A)$,
\item"$(\circ)$"
the sets
$B_i:=T^iA_i$ are mutually disjoint Borel subsets of $B$,
\item"$(\circ)$"
$|\rho_\mu(x,T^ix)-w|<\delta/2$ for  a very little
$0<\delta<\epsilon$ at all $x\in A_i$  and
\item"$(\circ)$"
$1/2\mu(B)<\sum_{|i|<N}\mu(B_i)<3/4\mu(B)$.
\endroster
Now find mutually
disjoint clopen subsets $A_i^0\subset A$ such that
\roster
\item"$(\circ)$"
$\mu(A_i^0\triangle A_i)<\delta\mu(A_i)$,
\item"$(\circ)$"
the subsets
$B_i^0:=T^iA_i^0$, $|i|<N$, are mutually disjoint (clopen) subsets of
$B$ and
\item"$(\circ)$"
$|\rho_\mu(x,T^ix)-w|<\delta$ for all $x\in A_i^0$.
\endroster
We now
define a homeomorphism---a `piece' of $S$---from a clopen subset
$A_0:=\bigsqcup_{|i|<N}A_0^i$ of $A$ onto a clopen subset
$B_0:=\bigsqcup_{|i|<N}B_0^i$ by setting $Sx:=T^ix$ if $x\in A_i^0$.
Notice that
\roster
\item"$(\circ)$"
$1/2<\mu(A_0)/\mu(A)<3/4$, $1/2<\mu(B_0)/\mu(B)<3/4$,
\item"$(\circ)$"
$|\rho_\mu(x,Sx)-w|<\delta<\epsilon$ for all $x\in A^0$, and
\item"$(\circ)$"
 the map
$A_0\ni x\mapsto\rho_\mu(x,Sx)$ is continuous.
\endroster
This completes the
first step of the inductive procedure. Now repeat the same for the
pair of clopen subsets $A\setminus A_0$ and $B\setminus B_0$ (notice
that the real $\log(\mu(B\setminus B_0)/\mu(A\setminus A_0))$ belongs
to $r(\rho_\mu)$ and it is very close to $w$) and so on. It remains
to concatenate the pieces of $S$ that were constructed on all steps
to obtain a homeomorphism $S$ from an open subset $A'$ of $A$,
$\mu(A')=\mu(A)$ to an open subset $B'$ of $B$, $\mu(B')=\mu(B)$. \qed
\enddemo

The following statement is a $III_1$-analogue of Lemma~2.4. We omit
proof of it since it coincides almost verbally with that of Lemma~2.4
(just replace the reference to Lemma~3.1 with the reference to
Lemma~2.3).

\proclaim{Lemma 3.2} Let $(\mu, T)$ be as in Lemma~3.1. Let
$D,\widetilde
D$ be two countable discrete spaces and $B$ a clopen subset of $X$.
Let
$\nu$ be a finite non-degenerated measure on $I_r\times \widetilde
D$. Given $\epsilon>0$ and continuous onto maps $u:B\to D$ and
$\widetilde
v:I_r\times\widetilde D\to D$ such that $\mu\circ
u^{-1}=\nu\circ\widetilde v^{-1}$, there are an $r$-uniform
equivalence relation $\Cal S\subset\Cal R_T$ on $B$, an $\Cal
S$-fundamental subset $\widetilde B\subset B$, an $(\Cal S,\widetilde
B)$-splitting homeomorphism $v:I_r\times\widetilde B\to B$ and a
measure
preserving continuous onto map $\widetilde u:\widetilde
B\to\widetilde D$
such that the following diagram commutes
$$
\CD
 (B,\mu) @<{v}<<    (I_r\times\widetilde B,\mu\circ v)\\
 @V{u}VV    @VV{\text{id}\times\widetilde u}V\\
 (D,\nu\circ \widetilde v^{-1}) @<{\widetilde v}<<
(I_r\times\widetilde D,
\nu )
\endCD
$$
and $|\rho_\mu(v(j,\widetilde b),v(i,\widetilde
b))-\log(\nu(j,\widetilde u(\widetilde b))/\nu(i,\widetilde
u(\widetilde b)))|<\epsilon$ for all $i,j\in I_r$ and $\widetilde
b\in\widetilde B$.
Here  $\widetilde B$ and $\widetilde D$ are identified with the
`bottom' levels $\{0\}\times \widetilde B$ and $\{0\}\times
\widetilde D$ of $I_r\times \widetilde B$ and $I_r\times \widetilde
D$ respectively equipped with the induced measures.
\endproclaim

We now prove the main result of this section.

\proclaim{Theorem 3.3} Let $(X,\tau,\mu,T)$
and $(X',\tau',\mu',T')$ be two dynamical systems of type $III_1$.
They are almost continuously orbit equivalent.
\endproclaim

\demo{Sketch of the proof}
The structure of the proof is the same as in  the proof of
Theorem~2.5.
We use the same notation as there and emphasize only the new ideas in
the proof. Fix a sequence $\epsilon_n>0$ such that
$\sum_{n=1}^\infty\epsilon_n<\infty$. As in the proof of Theorem~2.5,
we will construct a commuting diagram \thetag{2-2}. However the
vertical arrows (and hence the orbit equivalence
$(\phi')^{-1}\circ\phi$) will no longer be measure preserving.
Nevertheless it will be nonsingular with continuous Radon-Nikodym
derivative.

We replace $(\bigtriangledown)$ with
\roster
\item"$(\bigtriangledown)'$"
the Radon-Nikodym derivative on ${\Cal R}_n$ is `up to $\epsilon_n$'
equivariant under $p_n$, i.e. given $d\in D_n$ and $i\in I_{r_n}$,
the map $p_n^{-1}(d)\ni b\mapsto\rho_\mu(h_n(0,b),h_n(i,b))$ is
constant
up to $\epsilon_n$
 \endroster
if  $n$ is odd and
\roster
\item"$(\bigtriangledown)''$"
the Radon-Nikodym derivative on $\widetilde{\Cal R}_n$ is up to
$\epsilon_n$ equivariant under $\pi_n$, i.e. given
$(i_1,\dots,i_n,d_n)\in I_{r_1}\times\cdots\times I_{r_n}\times D_n$,
the map
$$
p_n^{-1}(d_n)\ni
b_n\mapsto\sum_{j=1}^n\rho_\mu(h_j(0,b_j),h_j(i_j,b_j))
$$
is constant up to $\epsilon_n$, where $b_{j-1}:=h_{j}(i_j,b_j)$,
$j=n,n-1,\dots,2$,
\endroster
 if $n$ is even.
 We also replace references to Lemma~2.4 with references to
Lemma~3.2. The rest of the construction is same.

Let us verify that $\mu\circ\phi^{-1}\sim\mu'\circ{\phi'}^{-1}$ with
continuous Radon-Nikodym derivative.
First, we note that $(X_0,\mu)=\projlim_{n\to\infty}(X_n,\mu_n)$,
where
$X_n:=I_{r_1}\times\cdots\times I_{r_n}\times B_n$ and $\mu_n$ is the
image of $\mu$ under the canonical homeomorphism of $X_0$ with $X_n$
(see the upper line in \thetag{2-2}). On the other hand, we have
 a sequence of  measures $\nu_n$ on $I_{r_n}\times D_n$ `intertwined'
by the maps $g_n$.
They lead naturally  to an `inverse' sequence of probability spaces
$$
(Z_1,\eta_1)@<\text{id}\times g_2<<
(Z_2,\eta_2)@<\text{id}\times\text{id} \times g_3<<\cdots,
$$
where $Z_n:=I_{r_1}\times\cdots\times I_{r_n}\times D_n$ (see the
middle line in \thetag{2-2}).
Let $\eta$ denote the inverse limit of $(\eta_n)_n$. This is a
probability measure on $Z$. Notice that
$\phi=\projlim_{n\to\infty}\phi_n$, where
$\phi_n=\text{id}\times\cdots\times\text{id}\times p_n$ maps $X_n$
onto $Z_n$. We now verify  that the sequence of continuous functions
$$
a_n:Z\ni z\mapsto
a_n(z):=\log(d\mu_n\circ\phi_n^{-1}/{d\eta_n})(z_n)\in\Bbb R
$$
converges uniformly on $Z$ along odd $n\to\infty$ . For this, we
introduce some notation.
Given $i\in I_{r_n}$ and $d\in D_n$, we set
$\nu_n(i|d):=\nu_n(i,d)/\nu_n(0,d)$. Then for all $i_1\in
I_{r_1},\dots, i_n\in I_{r_n}, d_n\in D_n$, we have
$$
\align
\eta_n(i_1,\dots,i_n,d_n) &=\nu_1(i_1|d_1)\nu_2(i_2|d_2)\cdots
\nu_{n-1}(i_{n-1}|d_{n-1})\nu_n(i_n,d_n) \text{ and}\\
\mu_n\circ\phi_n^{-1}(i_1,\dots,i_n,d_n)
&=\int_{p_n^{-1}(d_n)}\prod_{j=1}^n
\exp(\rho_\mu(h_j(0,b_j),h_j(i_j,b_j)))\,d\mu(b_n),
\endalign
$$
where $d_{j-1}=g_j(i_{j},d_{j})$ and $b_{j-1}=h_j(i_{j},b_{j})$, for
$j=n,n-1,...,2$.
Since the map $\text{id}\times p_n$ is measure preserving, we deduce
from $(\bigtriangledown)'$ that
$$
\rho_\mu(h_n(0,b),h_n(i,b))=\log(\nu_n(i|d))\pm\epsilon_n\text{ for
all }
b\in p^{-1}_n(d),\ i\in I_{r_n}, d\in D_n,
\tag{3-1}
$$
if $n$ is odd.
On the other hand, $(\bigtriangledown)''$ yields
$$
\int_{A}\prod_{j=1}^n
\exp(\rho_\mu(h_j(0,b_j),h_j(i_j,b_j)))\,d\mu(b_n)=
\frac{\mu\circ\phi_n^{-1}
(i_1,\dots,i_n,d_n)\mu(A)(1\pm\epsilon_n)}{\mu(p_n^{-1}(d_n))}
$$
for every $(i_1,\dots,i_n,d_n)\in I_{r_1}\times\cdots\times
I_{r_n}\times D_n$ and Borel subset $A\subset p_n^{-1}(d_n)$ if $n$
is even. This fact and \thetag{3-1} imply that
 $|a_{n+1}(z)-a_n(z)|\le 2\epsilon_n$ for all $z\in Z_0$ if $n$ is
odd. Hence $\mu\circ\phi^{-1}\sim\eta$ with continuous Radon-Nikodym
derivative. In a similar way one (changing the parity of $n$) can
show that $\mu'\circ{\phi'}^{-1}\sim\eta$ with continuous
Radon-Nikodym derivative.
\qed

\enddemo

{}

\enddocument